\documentclass{amsproc}

\usepackage{amsmath}
\usepackage{amssymb}
\usepackage{amsthm}
\usepackage[all,2cell]{xy}

\newtheorem{thm}{Theorem}[section]
\newtheorem{lemma}[thm]{Lemma}

\newtheorem{prop}[thm]{Proposition}

\newtheorem{defn}[thm]{Definition}

\newtheorem{rem}[thm]{Remark}
\newtheorem{rmk}[thm]{Remark}

\def\G{{\mathbb G}}

\def\P{{\mathbb P}}
\def\Q{{\mathbb Q}}

\def\Z{{\mathbb Z}}

\def\cI{{\mathcal I}}

\def\cO{{\mathcal{O}}}
\def\cQ{{\mathcal{Q}}}

\def\Q{{\mathbb{Q}}}
\def\G{{\mathbb{G}}}

\def\operatorname#1{\mathop{\rm #1}\nolimits}

\def\Chow{\operatorname{Chow}}

\def\Pic{\operatorname{Pic}}

\def\deg{\operatorname{deg}}

\def\det{\operatorname{det}}

\begin{document}

\title[Non extendability theorem for Grassmannians]{An extension of Fujita's non
extendability theorem for Grassmannians}

\author{Roberto Mu\~noz}
\address{Departamento de Matem\'atica Aplicada, ESCET, Universidad
Rey Juan Carlos, 28933-M\'ostoles, Madrid, Spain}
%\curraddr{Department of Mathematics and Statistics,
%Case Western Reserve University, Cleveland, Ohio 43403}
\email{roberto.munoz@urjc.es}
\thanks{Partially supported by the Spanish government project MTM2006-04785 and by MIUR (PRIN project: Propriet\`a
geometriche delle variet\`a reali e complesse)}

\author{Gianluca Occhetta}
\address{Dipartimento di Matematica, Universit\`{a} di Trento, Via Sommarive 14, I-38050, Povo (Trento), Italy
}
\email{gianluca.occhetta@unitn.it}
%\thanks{Partially supported by the Spanish government project MTM2006-04785.}

%    Information for third author
\author{Luis E. Sol\'a Conde}
\address{Departamento de Matem\'atica Aplicada, ESCET, Universidad
Rey Juan Carlos, 28933-M\'ostoles, Madrid, Spain}
\email{luis.sola@urjc.es}
%\thanks{Partially supported by the Spanish government project MTM2006-04785.}

%    General info
%\date{}


\begin{abstract}
In this paper we study smooth complex projective varieties $X$ containing a
Grassmannian of lines $\G(1,r)$ which appears as the zero locus of a section of a
rank two nef vector bundle $E$. Among other things we prove that the bundle $E$
cannot be ample.
\end{abstract}

\maketitle

\section{Introduction}\label{sec:intro}

\medskip

A natural problem in algebraic geometry is to study to which extent the geometry
of a smooth irreducible variety $X$ is determined by the geometry of its smooth
subvarieties $Y \subset X$, under certain positivity conditions on the embedding
$Y\subset X$. A typical result of this kind would
characterize $X$ (possibly saying that it cannot exist) by containing a
particular subvariety $Y \subset X$.

The classical setting in which the problem arose was the classification of
smooth projective embedded varieties $X \subset \P^N$ in terms of their smooth
linear sections $Y =X \cap \P^{N-k}$, for example the classification of low
degree embedded varieties. Later on, it evolved in different settings. For
instance one may impose positivity conditions on the
normal bundle $N_{Y/X}$.

If it is (generically) globally
generated then there exists a family of deformations of $Y$ sweeping out $X$.
That is the case, for example, of the
varieties swept out by linear subspaces of small codimension, see for
instance \cite{S} and \cite{NO}. In a recent paper (cf. \cite{MS}) the first and
third author have
dealt with embedded varieties $X\subset\P^N$ swept out
by codimension two Grassmannians, that may be regarded as a
projectively-embedded counterpart of this paper. The case of quadrics has been also studied, see  \cite{Fu} and
\cite{BIquadrics}.

If $N_{Y/X}$ is an ample line bundle, then it
is well known that, up to a birational transformation, $Y$
can be considered as an ample divisor on $X$. See \cite{HL} for a foundational reference
on ample subvarieties. With no assumption on the codimension, the hypothesis on $N_{Y/X}$ to be ample joint to some
topological assumptions constitute the setup of \cite{BdFL}. In that paper it is
shown how some structural maps (RC-fibrations, nef-value morphism, Mori contractions) of
$Y$ extend to $X$.

In the context of complex geometry, Lefschetz Theorem shows us how the topology
of $Y$ is reflected on the topology of a variety $X\subset\P^N$ containing $Y$
as a linear section. Moreover an extension of this result, due to Sommese (cf.
\cite{sommesefund} and \cite{sommese}), allows to work under weaker assumptions on the embedding
$Y\subset X$. In this way, Lefschetz-Sommese Theorem provides an important tool
for the type of problems we are considering here.

In this paper we will consider varieties $Y$ appearing as the zero locus of a
regular section of an ample vector bundle $E$ on $X$. An interesting survey on this matter has
been recently written by Beltrametti and Ionescu, see \cite{BI}. It deals mostly with the divisor case, but it also provides references for higher codimension. Among different results of
this kind, let us recall a theorem by T. Fujita (cf. \cite[Thm.~5.2]{F1}). It
states that, apart of the obvious cases, Grassmannians cannot appear as ample
divisors on a smooth variety. Our goal here is to show how this result can be
extended to codimension two:

\begin{thm}\label{thm:fujita} For $r \geq 4$ the Grassmannian of lines in $\P^r$
cannot appear as the zero locus of
a section of an ample vector bundle of rank two over a smooth complex projective
variety.
\end{thm}

Let us observe that the Grassmannian of lines in $\P^r$, say $\G(1,r)$, is
embedded naturally in $\G(1,r+1)$ as the zero locus of the universal quotient
bundle $\cQ$ which is not ample but globally generated. It is then natural to
look for a broader positivity assumption on $E$ in which this situation is
included. Taking in account Lefschetz-Sommese Theorem, it makes sense to
consider the notion of $k$-ampleness introduced by Sommese in
\cite[Def.~1.3]{sommese} (see also Definition \ref{def:k-ample} below). The main
result of this paper, from which Theorem \ref{thm:fujita} is a straightforward
corollary, is the following:

\begin{thm}\label{thm:main} Let $X$ be a smooth complex projective variety of
dimension $2r$ and
$G\subset X$ a subvariety isomorphic to the Grassmannian of lines in $\P^r$,
$r\geq 4$. We further assume that $G$ equals the zero set of a section of a
$(2r-4)$-ample vector bundle $E$ on $X$ of rank two. Then $X$ is isomorphic to
the Grassmannian of lines in $\P^{r+1}$ and $E$ is the universal quotient bundle
of this Grassmannian.
\end{thm}

Our proof relies on proving that the normal bundle of $G$ in $X$ must be
uniform, and on the classification of uniform vector bundles of low rank on
Grassmannians.

The structure of the paper is the following. In Section \ref{sec:preliminars} we
recall some generalities on Grassmannians, positive vector bundles and vanishing
results that we will use along the paper. In particular we find a lower bound on
the degree of $E$ in terms of the index of $X$. Moreover one may show that this
index is at most $\dim X-2$, a fact that is crucial in our argumentation; this
is the purpose of Section \ref{sec:ruleout}.  Section \ref{sec:uniform} deals
with the classification of uniform vector bundles on Grassmannians, and in
Section \ref{sec:E} we determine the possible values of the restriction $E|_G$.
In Section \ref{sec:proof} we present the proof of Theorem \ref{thm:main}, and
finally in Section \ref{sec:applications} we use the results in \cite{BdFL} in
order to derive from Theorem \ref{thm:fujita} a non-extendability result for
Grassmannian fibrations.

\noindent{\bf Acknowledgements:} We would like to thank Tommaso de Fernex for
his useful comments regarding Grassmannian fibrations.

\subsection{Conventions and definitions.}\label{ssec:conventions}
Along this paper $X$ will denote a smooth complex projective variety of
dimension $2r$ and
$G\subset X$ a subvariety isomorphic to the Grassmannian of lines in $\P^r$,
$\G(1,r)$, with $r\geq 4$. We further assume that $G$ equals the zero set of a
section of a  vector bundle $E$ on $X$ of rank two which is $(2r-4)$-ample in
the sense of Sommese, see Definition \ref{def:k-ample}. Denoting by $\cO(1)$ the
ample generator of $\Pic(G)\cong\Z$, the determinant of $E|_G$ is isomorphic to
$\cO(c)$, for some $c\in \Z$. We call $c$ the {\it degree of} $E$. The notation
$\cO(1)$ will be also used to denote the ample generator of a variety of Picard
number one and the tautological line bundle on a projective bundle. Subscripts
will be used if necessary.

Finally, on a Fano variety of Picard number one, a rational curve of degree one
with respect to $\cO(1)$ will be called a {\it line}. By definition, the family
of lines in $X$ is {\it unsplit}, i.e. the subscheme of $\Chow(X)$ parametrizing
them is proper. That amounts to say that a line is not algebraically equivalent
to a reducible cycle.

We will write $\mathbb{Q}^n$ (or just $\mathbb{Q}$ when its dimension is not
relevant) for a $n$-dimensional smooth quadric.

\section{Preliminaries}\label{sec:preliminars}
\subsection{Generalities on Grassmannians}\label{ssec:grass}
Let us recall some well-known facts on Grassmannians. We follow the conventions
of \cite{A}. As said before, the Grassmannian of lines in $\P^r$ is denoted by
$\G(1,r)$. We will denote by $\cQ$ the rank two {\it universal quotient bundle}
and by $\mathcal{S}^\vee$ the rank $r-1$ {\it universal subbundle}, related in
the universal exact sequence:
$$0 \to \mathcal{S}^\vee \to \cO^{\oplus r+1} \to \cQ \to 0.$$

The projectivization of $\cQ$ provides the universal family of lines in $\P^r$:
$$\xymatrix{&\P(\cQ)\ar[ld]_{p_1}\ar[rd]^{p_2}&\\ \G(1,r)&&\P^r.} $$
From right to left, this diagram may be thought of as the universal family of
$\P^{r-1}$'s in $\G(1,r)$. These $\P^{r-1}$'s have degree one with respect to
the Pl\"ucker polarization and their normal bundles in $\G(1,r)$ are isomorphic
to $T_{\P^{r-1}}(-1)$.
Finally we recall that the Chow ring of $\G(1,r)$ is generated by a well
determined type of cycles, called {\it Schubert cycles}. The generators in
dimension two are given by: the cycle parameterizing lines in a
$\P^3\subset\P^r$ passing by a point, and the cycle parameterizing lines in a
$\P^2\subset\P^r$ (we denote it by $\G(1,2)$). They are called $a$ and
$b$--planes, respectively.

\begin{rmk}\label{rmk:chern2}
{\rm In particular, the second Chern class of a vector bundle $E$ on $\G(1,r)$
is given by two integers, corresponding to the second Chern classes of the
restrictions of $E$ to the planes described above.}
\end{rmk}

\subsection{Positivity, topology and vanishing results}\label{ssec:topo}
The hypotheses on $X$ in \ref{ssec:conventions} impose severe restrictions on
its topology. In order to describe them explicitly let us recall the definition
of $k$-ampleness in the sense of Sommese, see \cite[Def.~1.3]{sommese}:

\begin{defn}\label{def:k-ample} {\rm Let $E$ be a semiample vector bundle over a
projective variety, i.e. $\cO_{\P(E)}(m)$ is free for $m$ big enough. The vector
bundle $E$  is said {\it $k$-ample} if every fiber of the morphism $\phi:\P(E)
\to \P(H^0(\P(E),\cO_{\P(E)}(m)))$ has dimension less than or equal to $k$.}
\end{defn}

In particular any $k$-ample vector bundle is nef and is ample if it is
$0$-ample.
Sommese's extension of Lefschetz Hyperplane Section Theorem  \cite[II,~Thm.~
7.1.1]{lazarsfeld} admits an extension to $k$-ample vector bundles, see
\cite[Prop.~1.16]{sommese} quoted in \cite[II,~Rmk.~ 7.1.9]{lazarsfeld}, which
applies to our case giving the following relations between the topologies of $X$
and $G$.

\begin{lemma}\label{lemma:lefschetz} Let $X$, $G$ and $E$ be as in
\ref{ssec:conventions}. The restriction map
$r:\Pic(X)\to\Pic(G)$ is an isomorphism.
\end{lemma}

\begin{proof} Denote by $r_i:H^{i}(X,\Z) \to H^{i}(G,\Z)$ the corresponding
restriction morphisms. By \cite[Prop.~1.16]{sommese} we get that $r_1$ is an
isomorphism and $r_2$ is injective with torsion free cokernel. Furthermore we
may compare the exponential sequences of $X$ and $G$ to get the following
diagram:
$$\xymatrix{ H^1(X,\Z)\ar[r] \ar[d]_{r_1} & H^1(X,\cO) \ar[r] \ar[d]_{r_{1,1}} &
\Pic(X) \ar[r] \ar[d]_{r}& H^2(X,\Z) \ar[d]_{r_2}\ar[r]& H^2(X,\cO)
\ar[d]_{r_{2,0}}
 \\
  H^1(G,\Z)\ar[r] & 0 \ar[r] & \Pic(G) \ar[r]^{\cong} & \Z \ar[r]& 0}$$
Since $r_1$ is an isomorphism and is compatible with the Hodge decomposition
then $r_{1,1}$ is an isomorphism. Since $r_2$ is injective and with torsion free
cokernel then it is an isomorphism and moreover $r_{2,0}$ is an isomorphism.
This implies that $r$ is an isomorphism.
\end{proof}

We will denote by $\cO_X(1)$ the ample generator of $\Pic(X)$, whose restriction
to $G$ is the Pl\"ucker line bundle. The degree of the canonical sheaf of $X$
equals $\deg(K_G)-\deg(E)$, that is $K_X=\cO(-r-1-c)$ and $X$ is Fano. Hence
Kobayashi-Ochiai Theorem \cite{KO} provides the bound
\begin{equation}\label{eq:c<=r} c\leq r. \end{equation}

Along this paper we will make use several times of the following variant of a
vanishing theorem due to Griffiths \cite[II,~Variant~ 7.3.2]{lazarsfeld}:

\begin{thm}\label{thm:vanishing}
 Let $M$ be a smooth complex projective variety of dimension $n$, $L$ an ample
line bundle on $M$ and $F$ a nef  vector bundle of rank $k$ on $X$, then:
$$H^i(M,\omega_M\otimes S^mF\otimes\det F\otimes L) =  0 \mbox{ for all }i>0,
m\geq 0. $$
\end{thm}

Applied to our setting, the previous theorem provides the following vanishing.

\begin{lemma}\label{lem:vanishing}
 Under the assumptions in \ref{ssec:conventions} and for every positive integer
$l$, it follows that
$$
 H^i(X,S^m E(l-r-1))=0, \mbox{ for all } i\geq 1, m\geq 0.
$$
\end{lemma}

Being $G$ the subscheme of zeroes of a section of the rank two vector bundle
$E$, the ideal sheaf of $G$ in $X$ has the following locally free presentation:

\begin{equation}\label{eq:presenta}
 0\rightarrow\det(E^\vee)\cong\cO(-c)\longrightarrow E^\vee\cong
E(-c)\longrightarrow\cI_{G/X}\rightarrow 0.
\end{equation}

Combining it with Lemma \ref{lem:vanishing} we immediately obtain:

\begin{lemma}\label{lem:h1(I)}
 With the assumptions of \ref{ssec:conventions}, the restriction maps
$$
H^0(X,\cO(k))\to H^0(G,\cO(k))
$$
are surjective for all $k>0$.
\end{lemma}
\begin{proof}
 In fact, it is enough to check that $H^1(X,\cI_{G/X}(k))=0$. Taking cohomology
on sequence (\ref{eq:presenta}), it suffices to show that
$H^1(X,E(k-c))=H^2(X,\cO(k-c))=0$. By Lemma \ref{lem:vanishing}, the first
vanishing holds whenever $k-c+r+1\geq 1$. Since $c\leq r$, see (\ref{eq:c<=r}),
that inequality is fulfilled for every positive $k$. For the second vanishing
note that since $\Pic(X) \cong \Z$, Kodaira vanishing implies that line bundles
on $X$ have no intermediate cohomology.
\end{proof}

Let us take a projective space of maximal dimension contained in $G$, say
$\P^{r-1} \cong M \subset G$, and denote by $E_M$ the restriction of $E$ to $M$.
Later on we will need to apply Theorem \ref{thm:vanishing} to $E_M$:

\begin{lemma}\label{lem:vanishing2}
 With the same assumptions as in \ref{ssec:conventions} and for every positive
integer $l$ it follows that:
$$
 H^i(M,S^m E_M(l+c-r))=0, \mbox{ for } i\geq 1, m\geq 0.
$$
\end{lemma}

%%%%%%%%%%%%%%%%%%%%%%%%%%%%%%%%%%%%%%%%%%%%%%%%%

\section{High index Fano varieties containing codimension two
Grassmannians}\label{sec:ruleout}

With the same assumptions as in \ref{ssec:conventions}, we will rule the cases
$c=r,\; r-1$ and $r-2$ out, which correspond to projective spaces, quadrics and
Del Pezzo varieties, respectively. In order to do that, it suffices to show that
$h^0(X,\cO(1))<r(r+1)/2=h^0(G,\cO(1))$ contradicting Lemma \ref{lem:h1(I)}.

In the case $c=r$ we get $h^0(X,\cO(1))=2r+1$ and hence it is smaller than
$r(r+1)/2$ whenever $r\geq 4$.

If $c=r-1$, $h^0(X,\cO(1))$ equals $2r+2$, which is smaller than $r(r+1)/2$ if
$r\geq 5$. The case $r=4$ would correspond to a smooth quadric $\Q^8\subset\P^9$
containing a Grassmannian $G\cong \G(1,4)$, embedded in $\P^9$ via the Pl\"ucker
map. But quadrics containing $\G(1,4)$ are given by $4\times 4$ pfaffians, hence
singular.

In order to rule out the case of Del Pezzo varieties, we will make use of
Fujita's classification (cf. \cite[8.11, p. 72]{F2}, see also \cite[V,1.12]{K}).
Being $\dim X=2r\geq 8$, the only possible values of $h^0(X,\cO(1))$ are
$2r,2r+1,2r+2$ and $2r+3$, which are smaller that $r(r+1)/2$ except in the
following cases:
\begin{itemize}
  \item $X$ is a smooth cubic hypersurface in $\P^9$ containing a Pl\"ucker
embedded Grassmannian $G\cong \G(1,4)$. Recall the notation on Grassmannians
established in \ref{ssec:conventions} and note that the normal bundle of $G$ in
$\P^9$ is $\cO(1)\otimes\wedge^2\mathcal{S}\cong\cO(2)\otimes\mathcal{S}^\vee$
where $\mathcal{S}^\vee$ denotes the universal subbundle, see for instance
\cite[Prop.~4.5.1]{manivel}. In particular, denoting by $E_G$ the restriction of
$E$ to $G$ we get the following exact sequence: $$0 \to E_G \to
\mathcal{S}^\vee(2) \to \cO(3) \to 0.$$ Tensoring by $\cO(-3)$ we get $H^1(G,
E_G(-3))\ne 0$. Now use Serre duality to get $h^1(G,E_G(-3))=h^5(G,E_G^\vee(3)
\otimes \omega_G)=h^5(G, E_G(1) \otimes \omega_G)$. Since $E_G(1)$ is ample we
get a contradiction with Le Potier Vanishing Theorem
\cite[Thm.~7.3.5]{lazarsfeld}.
 \item $X$ is a smooth complete intersection  of two quadrics $\mathbb{Q}_1$ and
$\mathbb{Q}_2$ in $\P^{10}$ containing a Pl\"ucker embedded Grassmannian $G\cong
\G(1,4)$. We may argue as before: observe on one hand that as a consequence of
Theorem \ref{thm:vanishing} we get that $h^1(G,E_G(-2))=0$. But on the other
hand taking cohomology on the following exact sequences we get the contradiction
$h^1(G,E_G(-2)) \ne0$:
     $$\begin{array}{c}\vspace{0.2cm}\xymatrix{
     0\ar[r] & E_G \ar[r] & N_{G/\P^{10}} \ar[r]& \cO(2)^{\oplus 2} \ar[r] &0
,\\
     }\\\xymatrix{0\ar[r] & S^\vee(2) \ar[r] & N_{G/\P^{10}} \ar[r]& \cO(1)
\ar[r] &0
      .}\end{array}
$$

    %In this case $G$ would lie in a hyperplane $H\subset\P^{10}$, and we get to
%a contradiction just by noting that quadrics in $\P^9$ containing $G$ have rank
%$6$, whereas the intersection $\mathbb{Q}_1\cap H$ has rank at least $9$.
\end{itemize}
As a corollary of what we have proved and recalling that $E$ is nef we get:
\begin{lemma}\label{lem:0<c<r-2} Under the assumptions of \ref{ssec:conventions}
we get that $0 \leq c <r-2$.\end{lemma}

%%%%%%%%%%%%%%%%%%%%%%%%%%%%%%%%%%%%%%%%%%%%%%%%%

\section{Uniform vector bundles on Grassmannianns}\label{sec:uniform}

Uniform vector bundles of low rank on Grassmannians have been classified by
Guyot, cf. \cite{G}. For the sake of completeness we present here a proof for rank
two vector bundles $E$ on $\G(1,r)$, using minimal
sections of $E$ over its lines. Although
we need only the case $r \geq 4$ we include a proof working for any $r \geq 2$.

Let us recall that a rank $k$ vector bundle $E$ on
$\G(1,r)$ is {\it uniform of type} $(a_1, \dots, a_k)$ ($a_1 \leq \dots\leq
a_k)$ if for any line $\ell \subset \G(1,r)$ the restriction of $E$ to $\ell$
splits as $\cO(a_1) \oplus \dots \oplus \cO(a_k)$. The result is the following:

\begin{prop}\label{prop:class} Every uniform rank two vector bundle $E$ on
$G:=\G(1,r)$ of
type $(0,1)$ is isomorphic either to
$\cO\oplus\cO(1)$ or to the universal bundle $\cQ$.
\end{prop}

\begin{proof} Note that for $r=2$ the result is due to Van de Ven (cf. \cite{VV},
\cite[Thm.~2.2.2]{OSS}), and we may assume
that $r\geq 3$.

%Consider the projective bundle $\pi:\P(E) \to G$ and denote by $V$ the family
% of lines in $G$. Each curve $\ell$ of $V$ admits a unique lifting to a curve
%in
% $\P(E)$ of $\cO(1)$ degree $0$ determined by the unique surjective morphism
% $E|_\ell\to\cO_\ell$. In this way $V$ defines an unsplit family of curves in
% $\P(E)$, that we denote also by $V$. Since $\Loc(V)\subset\P(E)$ dominates
%$G$,
% hence it is either an irreducible divisor or the whole $\P(E)$.
%
%If $\Loc(V)$ is a divisor then we claim that it provides a section of the
% projective bundle corresponding to a surjective morphism $E \to \cO(d)$. This
% implies that $E$ splits as a sum of line bundles.  Let us prove the claim.
%Note that $V_x$ is irreducible for every $x\in X$. For every $y\in \Loc(V)$ we
%get:
%$$
%\dim\Loc(V)_y\geq -K_{\P(E)}\cdot V-1+\codim(\Loc(V)\subset\P(E))=r.
%$$
%On the other side, since $\Loc(V)_{\pi^{-1}(\pi(y))}=\Loc(V)_{\pi(y)}$ is
% irreducible and $r$-dimensional, necessarily
% $\Loc(V)_{\pi^{-1}(\pi(y))}=\Loc(V)_y$, hence every curve of $V$ meeting the
% fiber $\pi^{-1}(\pi(y))$, does it at $y$. Hence we may assume that
% $\Loc(V)=\P(E)$.

First we show that there exists a family of linear subspaces of $G$ of maximal
dimension
verifying that $E|_G\cong\cO\oplus\cO(1)$. In fact, if $r\geq 4$, the
restriction of $E$ to a $\P^{r-1}$ is
isomorphic to $\cO\oplus\cO(1)$ by the classification of uniform vector bundles
on projective spaces (cf. \cite{EHS}, \cite[Thm.~3.2.3]{OSS}). For the case $r=3$ recall that the
Grassmannian $\G(1,3)$
contains two families of $\P^2$'s that we call $a$ and $b$-planes, see Section
\ref{ssec:grass}. Let us prove that the restriction of $E$ could not be
isomorphic to $T_{\P^2}(-1)$ for both families. If this occurs then $c_2(E)$
equals the union of two planes, one of each family, see Remark \ref{rmk:chern2}.
Assume by contradiction that this is the case. Consider two $a$-planes $a_1$ and
$a_2$ and
denote by $P$ their intersection and by $r$ the corresponding line
in $\P^3$. For every plane $M$ containing $r$ (determining a $b$-plane
$b_M$ containing $P$) we get two lines $r_1(M)=b_M\cap a_1$ and
$r_2(M)=b_M\cap a_2$. For each line $r_i$ ($i=1,2$) we get a lifting into
$\P(E)$ determined by the unique surjective map
$E|_{r_i(M)}\to\cO$. Denote by
$R_1(M)$ and $R_2(M)$ the intersections of these liftings with the fiber over
$P$. By hypothesis $E|_{a_i}\cong T(-1)$, hence the maps sending $M\to R_1(M)$
and $M\to R_2(M)$ are
isomorphisms from the set of planes containing $r$ to the fiber over
$P$. In particular there exists ${M_0}$ such that $R_1(M_0)=R_2(M_0)$. Now we
consider
the $b$-plane $b_{M_0}$. It contains two lines whose distinguished liftings
meet at one point. Then the restriction of $E$ to $b_{M_0}$ cannot be
$T(-1)$.

Recall that the family of $\P^{r-1}$'s of the previous paragraph is
parameterized by a projective space $\mathcal{M} \cong \P^r$. Each element of
this family admits a lifting to $\P(E)$ given by the unique surjective morphism
$E|_{\P^{r-1}}\to\cO$ and we have the following diagram:
$$\xymatrix{&\P(\cQ)\ar[ld]_{p_2}\ar[rd]^{p_1}\ar[r]^g&\P(E)\ar[d]^{\pi}\\
\mathcal{M} \cong \P^r& &\G(1,r)}$$
where $\cQ$ stands for the universal quotient bundle on $\G(1,r)$.

Now consider the restriction of $E$ to any $\G(1,2) \subset G$. The restriction
$E|_{\G(1,2)}$ is either decomposable or isomorphic to $T(-1)$ by Van de Ven's
result. We claim that in the former case $E$ is decomposable. In fact take a
point $x\in G$ and two $\P^{r-1}$'s, say $M_1$ and $M_2$, passing by $x$. We may
find a $\G(1,2)$ meeting $M_1$ and $M_2$ in two lines. The (unique) lifting of
this two lines to
$\P(E)$ as curves of degree $0$ with respect to $\cO(1)$ meet in one point,
since the two lines lie in $\G(1,2)$ and $E|_{\G(1,2)}=\cO\oplus \cO(1)$. In
particular $g(\P(\cQ))$ meets the fiber $\pi^{-1}(x)$ in one point, hence
$\pi:\P(E)\to G$ has a section and so $E$ splits as a sum of line bundles.

From now on we assume that $E|_{\G(1,2)}\cong T_{\P^2}(-1)$ for any
$\G(1,2)\subset G$. Arguing as in the previous paragraph, we may prove that in
this case the map $g$ is surjective.
Moreover there cannot be two liftings of $\P^{r-1}$'s passing by the same
point of $\P(E)$. In fact, if this occurs, we push it down to $G$ and we find a
$\G(1,2)$ meeting the two
$\P^{r-1}$'s in two lines. But the (unique) lifting of this two lines to
$\P(E)$ as curves of degree $0$ with respect to $\cO(1)$ do not meet,
since they lie in $\G(1,2)$ and $E|_{\G(1,2)}\cong T_{\P^2}(-1)$, a
contradiction.

Summing up, the morphism $g:\P(\cQ)\to\P(E)$ is bijective, and the proof is
finished.
\end{proof}

%%%%%%%%%%%%%%%%%%%%%%%%%%%%%%%%%%%%%%%%%%%%%%%%%

\section{Determining $E$}\label{sec:E}

In this section we will prove the following:

\begin{prop}\label{prop:E} Under the assumptions of \ref{ssec:conventions} the
vector bundle $E$ verifies that the restriction of $E$ to $G$ is isomorphic to
$\cQ$, where $\cQ$ stands for the universal quotient bundle.
\end{prop}

We begin by studying the restriction $E_M$ of $E$ to a projective space of
dimension $r-1$, $\mathbb{P}^{r-1}\cong M \subset G$. As a consequence of the
upper bound on $c$ of \ref{lem:0<c<r-2} and of the numerical characterization of
rank two Fano bundles onto projective spaces, see for example \cite{APW}, we
get:

\begin{lemma}\label{lemma:splitting} Under the conditions above $E_M$ splits
either
as $E_M =\cO(1)^{\oplus 2}$ or as $E_M=\cO(1)\oplus \cO$.
\end{lemma}

\begin{proof} Take the projective bundle $\pi: \P(E_M) \to M$. Since
$-K_{\mathbb{P}(E_M)}=\cO(2) \otimes \pi^*\cO(r+1-c)$ then $\mathbb{P}(E_M)$ is
a Fano variety, i.e. $E_M$ is
a rank two Fano bundle.
Hence we can use the classification of rank two Fano bundles,
see \cite[Main Thm.]{APW} and \cite[Thm. (2.1)]{SW}, to get that either $E_M$
splits as a sum of line bundles or $r=4$, $c=2$ and $E_M=\mathcal{N}(1)$,
being $\mathcal{N}$ a null correlation bundle.
This last possibility is excluded by the bound $c<r-2$ of \ref{lem:0<c<r-2} so
that $E_M$ splits as
as a sum of line bundles.

Moreover the Bend and Break lemma leads to the following vanishing:
\begin{equation}\label{eq:vanishing}
H^0(M, E_M(-2))=0.
\end{equation}
In fact, consider the exact sequence:
$$
0\rightarrow T_{M}(-3)\cong N_{M/G}(-2)\longrightarrow
N_{M/X}(-2)\longrightarrow E_M(-2)\rightarrow 0.
$$
By Lemma \ref{lem:vanishing2} we get
$H^1(M,E_M(-2))=0.$
Taking cohomology in the Euler sequence tensored with $\cO(-2)$ we get that
$H^1(M,T_{M}(-3))=0$
and therefore $H^1(M,N_{M/X}(-2))=0$. In particular the subscheme
$\mathcal{M}_{\mathbb{Q}}$ of the Hilbert scheme parametrizing deformations of
$M$ in $X$ containing a fixed smooth quadric $\mathbb{Q}\subset M$ is
smooth at the point $[M]$ and its dimension equals $H^0(M,N_{M/X}(-2))$.

But $\mathcal{M}_{\mathbb{Q}}$ must be zero dimensional, otherwise given two
general points $p,q\in \mathbb{Q}$,
for every deformation $M_t$ of $M$ we could consider the line $\ell_t\subset
M_t$
joining $p$ and $q$.
Then a Bend and Break argument (cf. \cite[3.2]{De}) provides a reducible
 cycle $C$ algebraically equivalent to $\ell_t$, contradicting the fact that
$\ell_t$
has degree $1$ with respect to $\cO(1)$.
This implies that $H^0(M,N_{M/X}(-2))=0$, and so $H^0(M,E_M(-2))=0$, too.

Since $E_M$ is nef then, by the splitting of $E_M$ and (\ref{eq:vanishing}), we
get that $E_M=\cO(a_1)\oplus \cO(a_2)$ with
$0 \leq a_1 \leq a_2 \leq 1$ and $a_1+a_2=c$.
Hence $c\leq 2$, being $E_M=\cO(1)^{\oplus 2}$ if $c=2$ and $E_M=\cO(1)\oplus
\cO$
if $c=1$. If $c=0$ then denote by $E_G$ the restriction of $E$ to $G$.
The rank two vector bundle is uniform with respect to the family of lines and
in fact it is trivial, see \cite[(1.2)]{AW}. This contradicts the fact that the
Picard number of $X$ is one, see
\cite[Lemma 3.6]{MS}. \end{proof}

Now we can complete the proof of Proposition \ref{prop:E}.

\begin{proof}
First we prove that the case $c=2$ cannot occur. Recall that $E_G$ stands for
the restriction of $E$ to $G$.
As $E_M=\cO(1)^{\oplus 2}$ we get that for any line $\ell \subset G$ the
restriction of
$E_G$ to $\ell $ is $\cO_{\ell}(1)^{\oplus 2}$.
This implies uniformity of $E_G$ with respect to the family of lines and
moreover
$E_G=\cO(1)^{\oplus 2}$, \cite[(1.2)]{AW}. Consider the exact sequence of
(\ref{eq:presenta})
\begin{equation}
\label{eq:presntc=2}0 \to\cO(-2) \to E(-2) \to \mathcal{I}_{G/X} \to 0
\end{equation}
and tensor it by $E(-1)$ to get:
$$0 \to\ E(-3) \to E\otimes E(-3) \to E \otimes \mathcal{I}_{G/X}(-1) \to 0.$$
By the usual decomposition $E \otimes E \cong S^2 E \oplus \wedge^2 E$
and the vanishing of Lemma \ref{lem:vanishing} we get $h^1(X,E \otimes
\mathcal{I}_{G/X}(-1))=0$.
Now consider the exact sequence $$0 \to E \otimes \mathcal{I}_{G/X}(-1) \to
E(-1) \to E_G(-1)=
\cO^{\oplus 2} \to 0$$ to get that $h^0(X,E(-1))\geq 2$ and that $E(-1)$ is
generically globally generated. Hence, see \cite[Lemma 3.5]{MS},
$E(-1)=\cO^{\oplus 2}$. Tensoring the exact sequence (\ref{eq:presntc=2}) by
$\cO(1)$ we observe that $\mathcal{I}_{G/X}(1)$ is globally generated. This
implies, see \cite[Cor. 1.7.5]{BSbook}, that there exists a smooth element in
the linear system $|\cO(1)|$ containing $G$, which contradicts
\cite[Thm.~5.2]{F1}.

If $c=1$ we have shown in Proposition \ref{prop:class} that $E_G$ is either as
stated or splits as $E_G=\cO\oplus \cO(1)$. If $E_G$ splits, exactly as in the
proof of the case $c=2$, we get $H^0(X,E(-1))\ne 0$. But this is a
contradiction: in fact the exact sequence of (\ref{eq:presenta}) $$0 \to\cO(-1)
\to E(-1) \to \mathcal{I}_{G/X} \to 0$$ gives $H^0(X,E(-1))=0$. This concludes
that $E_G=\cQ$.
\end{proof}

\section{Proof of the main Theorem}\label{sec:proof}

Let us give the proof of Theorem \ref{thm:main}.

\begin{proof} Let us recall that as a consequence of what we proved in the
Section \ref{sec:E} we can suppose that $c=1$ and that $E_G=\cQ$.
Consider the projective bundle $\pi:\P(E)\to X$, which is a Fano variety. Recall
that $E$ is nef by hypothesis and not ample as $c=1$. Hence we get that for $m$
big enough the linear system $|\cO(m)|$ defines an extremal ray contraction
$\varphi$ leading to the following diagram:
$$
\xymatrix{\P(E)\ar[d]_{\pi}\ar[r]^{\varphi}&Z, \\X}
$$ where $Z$ is normal. For $G \subset X$ we get that $E_G=\cQ$ so that, taking
care of the Mori cone of $\P(E_G)$, the following diagram appears:
\begin{equation}\label{eq:diagramm}\xymatrix{ & \P^r \ar[dr]^{f}\\ \P(\cQ)
\ar@{^{(}->}[r]\ar[ur]^{\varphi_1}\ar[d]_{\pi_1}  & \P(E)
\ar[d]_{\pi}\ar[r]^{\varphi} & Z,\\ G \ar@{^{(}->}[r]& X}
\end{equation}
being $\pi_1$ and $\varphi_1$ the corresponding contractions of $\P(\cQ)$ and
$f$ finite onto its image, which implies that $\dim Z \geq r$.

Now we claim that the general fiber $F$ of $\varphi$ is isomorphic to $\P^{r}$.
In fact, $F$
is irreducible and smooth by Bertini's Theorem and, if it is not a single point,
adjunction formula tells us that
$$-K_{F}=-K_{\P(E)}|_F=\pi^*\cO(1)^{\otimes(r+1)}.$$
But $\pi|_{F}$ is finite, hence $\pi^*\cO(1)|_{F}$ is ample and the above
formula implies that, if not a point, $F$ is a Fano manifold of index greater
than or equal to $r+1$. Recall that $\dim F \leq r+1$, hence either $F$ is a
point or $F \cong \P^r$ and $\pi^*\cO(1)|_F=\cO(1)$ or $F$ is a smooth quadric
of dimension $r+1$. In order to exclude the first and the last possibility let
us introduce some notation. Since $N_{G/X}=E_G=\cQ$, which is globally
generated, then there exists a $(r+1)$-dimensional irreducible variety
$\mathcal{G}$ parameterizing deformations of $[G]$ which in fact contains the
point corresponding to $G$, say $[G] \in\mathcal{G}$, as a smooth point. The
family $\mathcal{G}$ dominates $X$. By rigidity of Grassmannians, the general
point $[G']\in \mathcal{G}$ is isomorphic to $\mathbb{G}(1,r)$. Moreover
$E_{G'}$ is nef and its Chern polynomial is that of $E_{G}$. In particular it is
uniform so that $E_{G'}=\cQ$, see Proposition \ref{prop:class}. Thus, for the
general point $y \in \P(E)$ there exists $[G_y] \in \mathcal{G}$ such that $y
\in \P(E_{G_y})$ and provides a diagram as the one of (\ref{eq:diagramm}).
Therefore \begin{equation}\label{eq:linearfiber} \varphi_1^{-1}f^{-1}(\varphi(y))
\supset \P^{r-1} \subset F \end{equation} and this inclusion excludes the
possibility of $F$ to be a point or a smooth quadric, being $r\geq 4$. Summing
up we have shown that $F\cong \P^r$ and $\pi^*\cO(1)|_F=\cO(1)$ so that
$\pi(F)\cong\P^r \subset X$. Moreover, since the fibers of $\varphi$ dominates
$X$ via $\pi$, then the normal bundle $N_{\pi(F)/X}$ is generically globally
generated.

We claim that $N_{\pi(F)/X}=T_{\P^r}(-1)$.  Consider the Euler sequence
$$ 0 \to \cO \to \pi^*(E^\vee) \otimes \cO(1) \to T_{\P(E)/X} \to 0$$ and
restrict it to $F$ to get that $$T_{\P(E)/X} \otimes \cO_F=\cO(-1).$$ Then,
since $\pi$ is an isomorphism, identifying isomorphic objects, we get the
following diagram:
$$\xymatrix{  &  & 0 \ar[d] & 0 \ar[d] &
 \\
  &  & T_{\P(E)/X} \otimes \cO_F \ar[r]^{\cong} \ar[d] & \cO(-1) \ar[d] &
 \\
0\ar[r] & T_F \ar[r] \ar[d]_{\cong} & T_{\P(E)} \otimes \cO_F \ar[r] \ar[d] &
N_{F/\P(E)}\cong\cO^{\oplus r+1} \ar[r] \ar[d] & 0
 \\
 0\ar[r]  &  T_{\pi(F)} \ar[r] & T_X \otimes \cO_{\pi(F)} \ar[r]\ar[d] &
N_{\pi(F)/X} \ar[r]\ar[d] & 0.
\\ &  & 0  & 0 &
 }$$ The last vertical sequence is that of Euler and $N_{\pi(F)/X}=T_{\P^r}(-1)$
as claimed.

Now we claim that $\varphi$ is equidimensional. Let us suppose the existence of
a fiber $F_0$ such that $\dim(\pi(F_0))>r$. Recall that for the general point $x
\in X$ there passes the image by $\pi$ of a general fiber $F$ of $\varphi$ and
moreover $\pi(F)\cong \P^{r}$, $N_{\pi(F)/X}=T_{\P^{r}}(-1)$ and
\begin{equation}\label{eq:selfintersection} c_r(N_{\pi(F)/X})=1.\end{equation}
Hence there exists a component $\mathcal{M}$ of the Hilbert scheme of $\P^r$'s
in $X$ containing $[\pi(F)]$ as a smooth point and sweeping out $X$.
Through the general point $x \in \pi(F_0)$ there exists $[M] \in \mathcal{M}$
such that $x \in M \cong \P^r \subset X$. Since $\dim(M \cap \pi(F_0)) \geq 1$
and $E|_M=\cO\oplus \cO(1)$ then the intersection $M \cap \pi(F_0)$ admits a
unique section into $\P(E)$ contracted by $\varphi$. It follows that $F_0$
intersects the only section $M_0$ over $M$ contracted by $\varphi$ so that it
contains it, i.e. $M_0 \subset F_0$. Now consider a general $y \in \P(E)$.
Recall that $F_y=\varphi^{-1}(\varphi(y)) \cong \P^r$ and $[\pi(F_y)] \in
\mathcal{M}$. Moreover, since $E|_{\pi(F_y)}=\cO \oplus \cO(1)$ then $F_y$ is
the unique section of $E|_{\pi(F_y)}$ contracted by $\varphi$. But now observe
that as a consequence of the selfintersection formula and
(\ref{eq:selfintersection}) any element in $\mathcal{M}$ is meeting $M$ and
therefore $\pi(F_y) \cap M\ne \emptyset$ which in particular gives $\pi(F_y)
\cap \pi(F_0) \ne \emptyset$. But this leads to the contradiction $F_y  \subset
F_0$ .

From the fact that $\varphi$ is equidimensional it follows that $\varphi: \P(E)
\to Z$ is a $\P^r$-bundle, that is all fibers are linear and $\varphi$ is
providing the structure of projective bundle, see \cite[2.12]{F3} quoted in
\cite[Prop.~3.2.1]{BSbook}. In particular $Z$ is smooth.

Recall that $\varphi$ is defined by the system $|\cO(m)|$. We claim that we may
assume $m=1$. In fact, take  $x \in G\subset X$ and the fiber of $\pi$ over it,
that is $\ell_x \cong \P^1=\pi^{-1}(x)$. Consider $y \in \ell_x$ and the fiber
$F_y$ of $\varphi$ through $y$. Now observe that $[\pi(F_y)] \in \mathcal{M}$
and that $F_y$ corresponds to the only section of $E|_{\pi(F_y)}^\vee$. Then
$F_y \cap \ell_x=\{y\}$ so that $\varphi|_{\ell_x}$ is a one-to-one map from
$\P^1$ onto its image in $Z$. Hence the restriction of $f$ to
$\varphi_1(\ell_x)$ is an isomorphism, for every $x\in G$. Since $G$
parametrizes all the lines of $\P^r$, it follows that $f$ itself is an
isomorphism. Therefore we may consider $\P^r$ as an effective divisor in the
smooth variety $Z$. Since $\Pic(Z)=\Z$, then $\P^r\subset Z$ is ample and
Kobayashi-Ochiai Theorem tells us that $Z\cong\P^{r+1}$ and
$\cO_Z(\P^r)\cong\cO_{\P^{r+1}}(1)$. In particular, fibers of $\pi$ map onto
lines of $Z\cong\P^{r+1}$.

The next step in the proof is to observe that through any two points $x, y\in
X$ there cannot pass two elements of $\mathcal{M}$. In fact, by the self
intersection formula and (\ref{eq:selfintersection}) it holds that two possible
different elements $M_1, M_2$ of $\mathcal{M}$ through $x$ and $y$ must meet in
a positive dimensional subvariety $P=M_1 \cap M_2$. But $E_{M_i}=\cO \oplus
\cO(1)$ for $i=1,2$ so that, exactly as in the proof of the equidimensionality
of $\varphi$, the corresponding unique sections $\P^r \cong F_i \subset \P(E)$
such that $\pi(F_i)=M_i$ are going to the same point by $\varphi$, contradicting
the fact that $\varphi: \P(E) \to Z$ is a $\P^r$-bundle.

Recall that $\ell_x:=\varphi(\pi^{-1}(x))$ is a line in $Z\cong\P^{r+1}$ for all
$x\in X$. This provides a map $g: X \to \mathbb{G}(1,r+1)$ sending $x$ to
$\ell_x$. Since $X$ and $\mathbb{G}(1,r+1)$ are smooth of Picard number one then
 we conclude the proof of the theorem by showing that $g$ is surjective and
generically injective. It is then enough to prove that for the general $x \in X$
there is no $y \in X$ different from $x$ such that $\ell_x= \ell_y$. Suppose on
the contrary the existence of such $y \in  X\setminus\{x\}$. For any point $z
\in \ell_x$ we get that $\varphi^{-1}(z)=\P^{r}$ is meeting the lines
$\pi^{-1}(x)$ and $\pi^{-1}(y)$. This implies that $\pi(\varphi^{-1}(z))$ is the
only element $\P^r =M \in \mathcal{M}$ through $x$ and $y$. This provides a one
dimensional family of sections of $E_M^\vee$, which is a contradiction.
\end{proof}

\begin{rmk}\label{rmk:picardhypothesis} {\rm Let us remark that, as has been
seen in the course of the proof, the hypothesis on the $(2r-4)$-ampleness of $E$
can be substituted by the hypothesis on the restriction map
$r:\Pic(X)\to\Pic(G)$ to be an isomorphism. Note that $G$ appears as the zero
set of a $(2r-2)$-ample vector bundle on, for instance, the product
$X=G\times\P^2$, but
$\Pic(X)\neq\Z$. A similar situation appears by considering the
desingularization of a cone over $G$ with vertex a line. We do not know yet of
any example in which $E$ is $(2r-3)$-ample and the restriction $r$ is not an
isomorphism.}
\end{rmk}

\section{Low values of $r$}

The case $r=3$ can be seen as a particular case of the general problem of
quadrics appearing as the zero locus of sections of
positive rank two vector bundles. This is well understood in the case in which
$E$ is ample \cite{LM} (in fact in any codimension). Here we can prove the
following:

\begin{prop}\label{prop:quadrics} Let $X$ be as smooth complex projective
variety of dimension
$n \geq 6$. Suppose the
existence of a rank two nef vector bundle $E$ on $X$ and a section of $E$ vanishing on a
smooth quadric
$\mathbb{Q} \subset X$.
If the restriction map $r:\Pic(X) \to \Pic(\mathbb{Q})$ is an isomorphism then
$(X,E)$ is either
\begin{itemize}
\item $(\P^n, \cO(2) \oplus \cO(1))$, or

\item $(\mathbb{Q}, \cO(1) \oplus \cO(1))$, or

\item $(\mathbb{G}(1,4), \cQ)$.
\end{itemize}
\end{prop}

\begin{proof} Denote as usual by $\cO(1)$ the ample generator of $\Pic(X)$ and
by $c$ the degree
of the determinant of $E$. Recall that since $E$ is nef and has a section
vanishing on $E$ then $c>0$.
Now use adjunction formula to get that $K_X=\cO(-(n-2)-c)$. This implies that
either
$c=3$ and $X=\P^n$ or $c=2$ and $X=\mathbb{Q}$ or $c=1$. Hence we can suppose
that $c=1$
which means that $X$ is a Del Pezzo Variety. Now we apply \cite[Prop.~4.5]{MS}
to get that $X$
is $\mathbb{G}(1,4)$. If $X \cong \mathbb{G}(1,4)$ then $E$ either splits as a
sum of line bundles or $E \cong \cQ$, see Proposition \ref{prop:class}. But in
case $E=\cO \oplus \cO(1)$ there are no sections vanishing on a codimension two
variety and the result follows.
\end{proof}

The case $r=2$ can be seen as a particular case of the general problem of linear
spaces appearing as the zero locus of sections
of positive rank two vector bundles. See \cite{LM} for the case in
which $E$ is ample. Here we can prove the following:

\begin{prop}\label{prop:linear}
Let $X$ be a smooth complex projective variety of dimension $n \geq 4$. Suppose
the
existence of a rank two nef vector bundle $E$ on $X$ and a section of $E$ vanishing on a
linear space
$\P^{n-2} \subset X$.
If the restriction map $r:\Pic(X) \to \Pic(G)$ is an isomorphism then  $(X,E)$
is either
\begin{itemize}
\item $(\P^n, \cO(1) \oplus \cO(1))$,

\item $(\mathbb{G}(1,3), \cQ)$.
\end{itemize}
\end{prop}

\begin{proof} With the same notation as before we get by adjunction that
$K_X=\cO(-(n-1)-c)$. Then
either $c=2$ and $X \cong \P^n$ or $c=1$ and $X$ is a smooth quadric so that
$\dim(X) \leq 4$, in
fact equal by hypothesis.
Since $E$ is uniform then either $E=\cO \oplus \cO(1)$ and no section vanishes
in a codimension two subvariety or
$E=\cQ$ and we conclude.
\end{proof}

\begin{rmk}\label{rmk:changeofhypothesis} {\rm For Propositions
\ref{prop:quadrics} and \ref{prop:linear} let us remark that if we impose on $E$
to be $(n-4)$-ample then we get that the restriction morphism $r:\Pic(X) \to
\Pic(G)$ is an isomorphism, exactly as in Lemma \ref{lemma:lefschetz}.}
\end{rmk}

\section{Fibrations in Grassmannians of lines}\label{sec:applications}

Inspired by \cite[Def.~5.1]{BdFL} we can give the following definition.

\begin{defn}\label{def:Grassmannianfibration} {\rm A surjective morphism $\pi:Y
\to Z$ between
 a smooth projective variety $Y$ and a normal projective variety $Z$ is called a
${\mathbb G}(1,r)$-{\it fibration} if $\pi$ is an elementary Mori contraction
and there is a line bundle
$L$ on $Y$ such that the general fiber $G$ of $\pi$ is isomorphic to ${\mathbb
G}(1,r)$ and $L|_G$
is the Pl\"ucker line bundle.}
\end{defn}

If $\dim(Z) \leq 2$ it suffices to check this hypothesis on a fiber:

\begin{lemma}\label{lem:codimensiontwofibrations} Let $\pi:Y \to Z$ be a
morphism between a smooth projective variety $Y$ and a normal projective variety
$Z$ such that
$\dim(Z) \leq 2$. If there exists $z$ a smooth point of $Z$ and $L \in \Pic(Y)$
such that $G:=\pi^{-1}(z)$ is
isomorphic to $\mathbb{G}(1,r)$ and $L|_{G}$ is the Pl\"ucker line bundle then
$\pi:Y \to Z$ is a
$\mathbb{G}(1,r)$-fibration. Moreover $Z$ is smooth and all smooth fibers of
$\pi$ are isomorphic to $\mathbb{G}(1,r)$.
\end{lemma}

\begin{proof} Since the normal bundle $N_{G/Y}$ is trivial then, in particular
is generically globally generated and its determinant is also trivial.
Up to replacing $L$ with $L\otimes\pi^*A$, with a suitable ample line bundle $A$
on $Z$ we may assume that $L$ is ample and we can apply \cite[Lemma 2.5]{MS} to
get that $\pi$ is the contraction of an extremal ray. Moreover, by \cite[Cor.
1.4]{AW}, $\pi$ is equidimensional  and $Z$ smooth. By rigidity of
Grassmannians, see for instance \cite{HM},
any smooth fiber is isomorphic to $\mathbb{G}(1,r)$ and the lemma follows.
\end{proof}

\begin{prop}\label{prop:fibrations} For $r \geq 4$ a $\mathbb{G}(1,r)$-fibration
$Y$ cannot appear either as an ample divisor
or as the zero locus of a section of a rank two ample vector bundle $E$ over a
smoooh projective variety $X$.
\end{prop}

\begin{proof} Suppose on the contrary that $Y \subset X$ appears as the zero
locus of
a section of $E$. Then, by Lefschetz-Sommesse Theorem, the restriction map from
$\Pic(X)$ to $\Pic(Y)$ is an isomorphism. Hence we may use \cite[Thm 4.1]{BdFL}
to
get a diagram:
$$
\xymatrix{Y\ar[d]_{\pi}\ar@{^{(}->}[r] & X \ar[d]_{\phi} \\ Z \ar[r]^{\delta} &
S,}
$$
where $\phi$ is an elementary Mori contraction on
$X$ and $\delta$ is a finite morphism.

Consider a general point $s \in S$ and denote by $F_s=\phi^{-1}(s)$ the fiber of
$\phi$ over $s$, which is connected. Since $\delta$ is finite then
$\delta^{-1}(s)=\{z_1,\dots,z_d\}$. Denote by $G_i=\pi^{-1}(z_i)$, $1 \leq i
\leq d$, so that $F_s \cap Y=G_1 \cup \dots \cup G_d$, where $G_i \cap
G_j=\emptyset$ for $i \ne j$. Recall that $Y$ is defined as the zero locus of a
section of an ample vector bundle and, since $E|_{F_s}$ is ample then $F_s \cap
Y$ is also the zero locus of a section of an ample vector bundle. Then, by
Lefschetz-Sommesse Theorem, $F_s \cap Y$ is connected so that $d=1$, that is
$F_s \cap Y=G_1 \cong \mathbb{G}(1,r)$.
But  $\mathbb{G}(1,r)$ cannot appear
either as an ample divisor on $F_s$ by \cite[Thm.~5.2]{F1} or as the zero locus
of a
section of $E|_{F_s}$ by Theorem \ref{thm:fujita}. This concludes the result.
\end{proof}

\begin{rem}\label{rem:fibrations}
 {\rm Using \cite[Thm.~3.6]{BdFL}, a similar statement holds under different
hypotheses. We could have assumed that there exists an unsplit covering family
$V$ of rational curves in $X$ verifying the following: it restricts to a family
$V_Y$ covering $Y$ and the general equivalence class in $Y$ with respect to
$V_Y$ is isomorphic to $\G(1,r)$.}
\end{rem}

%%%%%%%%%%%%%%%%%%%%%%%%%%%%%%%%%%%%%%%%%%%%%%%%%%%%%%%%%%%%
\bibliographystyle{amsalpha}

\begin{thebibliography}{9999}

\bibitem[APW]{APW} Ancona, V., Peternell, T. and Wi\'sniewski. J.
  {\it Fano bundles and splitting theorems on projective spaces and quadrics},
Pacific J.  Math. {\bf 163},
  no. 1, 17-41 (1994).

\bibitem[AW]{AW} Andreatta, M. and Wi\'sniewski, J.
   {\it On manifolds whose tangent bundle contains an ample locally free
subsheaf}, Invent. Math. {\bf 146},
   209-217 (2001).

\bibitem[A]{A} Arrondo, E. {\it Subvarieties of Grassmannians}, Lecture Note
Series Dipartimento di Matematica Univ. Trento, 10 (1996).

\bibitem[BdFL]{BdFL} Beltrametti, M.C., de Fernex, T. and Lanteri, A. {\it Ample
subvarieties and
rationally connected fibrations}, Math. Ann. {\bf 341}, 897-926 (2008).

\bibitem[BI1]{BIquadrics} Beltrametti, M.C., and Ionescu, P.
{\it On manifolds swept out by high dimensional quadrics}. Math. Z. {\bf 260}, 229-234 (2008).

\bibitem[BI2]{BI} Beltrametti, M.C., and Ionescu, P. {\it A view on extending
morphisms from ample divisors}. To appear in Contemporary Mathematics.

\bibitem[BS]{BSbook} Beltrametti, M. C. and Sommese, A. J. {\it The
    Adjunction Theory of Complex Projective Varieties}, De Gruyter
  Expositions in Mathematics 16, De Gruyter, Berlin-New York, 1995.

\bibitem[BSW]{BSW} Beltrametti, M. C., Sommese, A. J. and
  Wi\'sniewski, J. {\it Results on varieties with many lines and their
    applications to adjunction theory}, Complex Algebraic Varieties
  (Bayreuth, 1990), Lecture Notes in Mathematics 1507, Springer,
  Berlin, 1992, pp. 16--38.

\bibitem[De]{De} Debarre, O.  {\it Higher-Dimensional Algebraic Geometry},
Springer-Verlag New York, 2001.

\bibitem[EHS]{EHS} Elencwajg, G., Hirschowitz, A. and Schneider, M.
{\it Les fibres uniformes de rang au plus {$n$} sur {${\bf P}\sb{n}({\bf C})$} sont ceux qu'on croit}.
Vector bundles and differential equations (Proc. Conf., Nice, 1979), Progr. Math., {\bf 7}, 37-63 (1980).

\bibitem[Fu]{Fu} Fu, B. {\it Inductive characterizations of hyperquadrics}. Math. Ann. {\bf 340},
no. 1, 185-194 (2008).

\bibitem[F1]{F1} Fujita, T. {\it Vector bundles on ample divisors}, J.
  Math. Soc. Japan {\bf 33}, no. 3, 405-414 (1981).

\bibitem[F2]{F3} Fujita, T. {\it On Polarized manifolds whose adjoint bundles are
not semipositive}.
in {\it Algebraic Geometry, Sendai 1985.} Adv. Stud. Pure Math. {\bf 10},
167-178
(1987).

\bibitem[F3]{F2} Fujita, T. {\it Classification Theories of Polarized
    Varieties}. London Mathematical Society Lecture Note Series, no.
  155. Cambridge University Press, Cambridge, 1990.

\bibitem[G]{G} Guyot, M. {\it Caract\'erisation par l'uniformit\'e des fibr\'es
universels sur la Grassmannienne}. Math. Ann. {\bf 270}, 47-62 (1985).

%\bibitem[Ha]{hartshorne} Hartshorne, R. {\it Algebraic geometry}.
%  Graduate Texts in Mathematics, No. 52. Springer-Verlag, New
%  York-Heidelberg, 1977.

\bibitem[H]{HL} Hartshorne, R. {\it Ample subvarieties of algebraic varieties}.
Lecture Notes in Mathematics, Vol. 156, Springer-Verlag, Berlin-New York 1970.

%\bibitem[Hw]{hwang-survey} Hwang, J-M. {\it Geometry f minimal
%    rational curves on Fano manifolds}, School on Vanishing Theorems
%  and Effective Results in Algebraic Geometry (Trieste, 2000), ICTP
%  Lect. Notes, vol. 6, Abdus Salam Int. Cent. Theoret. Phys., Trieste,
%  2001, pp. 335-393.

\bibitem[HM]{HM} Hwang, J-M. and Mok, N. {\it Rigidity of irreducible Hermitian
symmetric spaces
of the compact type under K\"ahler deformation}, Invent. Math. {\bf 131},
393-418 (1998).

%\bibitem[HM2]{HM2} Hwang, J-M and Mok, N. {\it Holomorphic maps from rational
%homogenous spaces of Picard number one onto projective manifolds}, Invent.
%Math. {\bf 136}, 209-231 (1999).

\bibitem[KO]{KO} Kobayashi, S., Ochiai, T. {\it Characterization of
    complex projective spaces and hyperquadrics}, J. Math. Kyoto Univ.
  {\bf 13}, 31-47 (1973).

\bibitem[K]{K} Koll\'ar J. {\it Rational curves on algebraic varieties},
Berlin, Springer-Verlag, 1996.

\bibitem[L]{lazarsfeld} Lazarsfeld, R. {\it Positivity in Algebraic
    Geometry II.} Springer-Verlag, Berlin-Heidelberg, 2004.

\bibitem[LM]{LM} Lanteri, A. and Maeda, H. {\it Ample vector bundles with
section vanishing on projective spaces or quadrics}, Int. J. Math.{\bf 6}, no.
4, 587-600 (1995).

\bibitem[M]{manivel} Manivel, L. {\it Gaussian maps and plethysm. Algebraic
geometry (Catania, 1993/Barcelona, 1994),} Lecture Notes in Pure and Appl.
Math., 200, Dekker, New York, 1998, 91-117.

\bibitem[MS]{MS} Mu\~noz, R. and Sol\'a-Conde, L. E. {\it Varieties swept out by
Grassmannians of lines}
to appear in Contemporary Mathematics.

\bibitem[NO]{NO} Novelli, C., Occhetta, G. {\it Projective
    manifolds containing a large linear subspace with nef normal
    bundle}, preprint 2008 arxiv: 0712.3406v2.

\bibitem[S]{S} Sato, E. {\it Projective manifolds swept out by large
dimensional linear spaces}, Tohoku Math. J. {\bf 49}, 299-321 (1997).

 \bibitem[OSS]{OSS}Okonek, C., Schneider, M. and Spindler, H. {\it
     Vector Bundles on Complex Projective Spaces}. Progress in
   Mathematics 3, Birkh\"auser, Boston, 1980.

\bibitem[So1]{sommesefund} Sommese, A. J. {\it On manifolds that cannot be ample
divisors}, Math. Ann. {\bf 221}, no. 1, 55--72  (1976).

\bibitem[So2]{sommese} Sommese, A. J. {\it Submanifolds of Abelian varieties},
Math.
Ann. {\bf 233}, no. 4 229-256 (1978).

\bibitem[SW]{SW} Szurek, M. and  Wi\'sniewski. J. {\it Fano bundles over $P^3$
and $Q_3$}, Pacific J. Math. {\bf 141}, no. 1 197-208 (1990).

\bibitem[VV]{VV}  Van de Ven, A. {\it On uniform vector bundles}. Math. Ann.
{\bf 195}, no. 4, 245-248 (1972).


%\bibitem[W]{Wis} Wi\'sniewski. J. {\it On Fano manifolds of large
%    index}, manuscripta math. {\bf 70}, 145-152 (1991).

\end{thebibliography}

\end{document}